\documentclass{article}
\pdfoutput=1

\usepackage{amsmath}
\usepackage{paralist}
\usepackage{graphics}
\usepackage{epsfig}
\usepackage{graphicx}
\usepackage{epstopdf}
\usepackage{float}
\usepackage{amssymb, bm}            
\usepackage{overcite}
\usepackage[nolist]{acronym}	    
\usepackage{dcolumn}	            
\usepackage{booktabs}
\usepackage{tabu}
\usepackage{subcaption}
\usepackage{cases}
\usepackage{etoolbox}
\usepackage{bm}
\apptocmd{\sloppy}{\hbadness 10000\relax}{}{}

\newcommand{\bbr}{\mathbb{R}}

\newcommand{\libzptrl}[2]{\dfrac{\partial#1}{\partial#2} }

\newcommand{\paren}[1]{\!\left(#1\right)\!}

\newcommand{\dint}{\displaystyle\int}

\begin{document}

\title{A Pseudospectral Approach to High Index DAE Optimal Control Problems}
	\maketitle
	
	\centerline{\scshape Harleigh C. Marsh,\footnote{Ph.D. candidate.}}
	\medskip
	{\footnotesize
		\centerline{Department of Applied Mathematics \& Statistics}
		\centerline{University of California, Santa Cruz}
		\centerline{ Santa Cruz, CA 95064, USA}
	} 
	
	\medskip
	
	\centerline{\scshape Mark Karpenko,\footnote{Research Associate Professor and corresponding author. e-mail: \tt{mkarpenk@nps.edu}}}
	\medskip
	{\footnotesize
		\centerline{Department of Mechanical and Aerospace Engineering}
		\centerline{Naval Postgraduate School}
		\centerline{Monterey, CA 93943, USA}
	}
	
	\medskip
	\centerline{and} 
		\medskip
	\centerline{\scshape Qi Gong\footnote{Associate Professor.}}
	\medskip
	{\footnotesize
		\centerline{Department of Applied Mathematics \& Statistics}
		\centerline{University of California, Santa Cruz}
		\centerline{ Santa Cruz, CA 95064, USA}
	}

	\vspace{1cm}

	\begin{abstract}
Historically, solving optimal control problems with high index differential algebraic equations (DAEs) has been considered extremely hard.  Computational experience with Runge-Kutta (RK) methods confirms the difficulties. High index DAE problems occur quite naturally in many practical engineering applications.  Over the last two decades, a vast number of real-world problems have been solved routinely using pseudospectral (PS) optimal control techniques. In view of this, we solve a ``provably hard,'' index-three problem  using the PS method implemented in DIDO$^\copyright$, a state-of-the-art MATLAB$^\circledR$ optimal control toolbox. In contrast to RK-type solution techniques, no laborious index-reduction process was used to generate the PS solution.  The PS solution is independently verified and validated using standard industry practices. It turns out that proper PS methods can indeed be used to ``directly'' solve high index DAE optimal control problems.  In view of this, it is proposed that a new theory of difficulty for DAEs be put forth.
	\end{abstract}

	\newpage
	\section{Introduction}

Systems of high dimensional nonlinear differential algebraic equations (DAEs) naturally arise when modeling physical processes in fields of mechanical, aerospace and electrical engineering\cite{campbell1997daes, kunkel2006differential, campbell2015flexibility}. Obtaining numerical solutions to high index DAE optimal control problems is widely considered to be difficult\cite{betts2013direct, campbell2016solving, campbell2016comments, campbell2007direct, engelsone2007direct}.  The computational experience for the difficulties largely stem from the use of Runge-Kutta (RK) methods for solving DAE optimal control problems.  One question we pose in this paper is: are DAE problems fundamentally difficult regardless of the computational method, or are the source of difficulties largely in the computational technique itself? The conventional wisdom points towards the former. In this paper we suggest that it might be the latter, and if so, a new DAE theory of difficulty is warranted.

An overview of our line of argument is as follows: Since the year 2007 when NASA implemented\cite{SIAMnews} a PS solution onboard the International Space Station, PS optimal control techniques have become the method of choice for solving ``NASA-hard'' problems\cite{bhatt:opm,TRACE-IEEE-Spectrum,TEI-JGCD-2011,Kepler-micro-slew,Karp-JWST}. As noted earlier, many of these practical problems are high-dimensional and differential-algebraic in nature.  DAE theory or index-reduction techniques have never been used to solve these problems. Although unlikely, it is quite possible that all these problems were fortuitously low-index DAE optimal control problems.  As a simple means to test the fundamental question related to the source of the hardship, a quick approach is to compare the numerically difficulties in solving a hard, high-index DAE problem using both RK and PS methods.  To this end, we investigate the high-index DAE optimal control considered by Campbell et al\cite{campbell2016solving,campbell2016comments}.  They showed that the suite of RK methods implemented in SOCX --- a state-of-the-art nonlinear-programming-based FORTRAN optimal control solver\cite{SOCX} --- failed to provide a solution to this problem without an index-reduction process\cite{campbell2016solving}. Similar difficulties were reported with other discretization methods and various other software tools\cite{campbell2016solving,campbell2016comments}.  Because DIDO$^\copyright$, a state-of-the-art MATLAB$^\circledR$ optimal control solver\cite{ross-book} was never used in their studies, this paper ``completes'' this effort.  As shown in this paper, DIDO is able to solve the challenge problem without any difficulty leading to the conclusion that a new theory of difficulty for DAEs might be warranted.

\section{The Challenge Problem}
The challenge problem posed by Campbell and Kunkel\cite{campbell2016solving} is to minimize the quadratic cost functional,
\begin{equation}\label{eqn:daeCostFun}
J(x, u)  = \dint_{0}^T cu^2 + d{(x- L\sin(t+\alpha))}^2 + d{(y- L\cos(t+\alpha))}^2\,dt
\end{equation}
subject to the constraints:
\begin{equation}\label{eqn:daeProblem}
(D):\begin{cases}
\ddot x = -\lambda x - a\dot x + uy\\
\ddot y = -\lambda y - a\dot y -g - ux\\
0 = x^2 + y^2 - L^2\\
x(0) = 0 \;\; y(0)=L \;\; \dot{x}(0)=0 \;\; \dot{y}(0)=0
\end{cases}
\end{equation}
The initial-value problem, $ (D) $, describes the motion of a simple pendulum of length $L$ in the Cartesian plane. The control, $ u $, is taken as the force applied tangentially to the motion of the pendulum bob. Nonlinear function, $ \lambda $, represents the internal force necessary to satisfy the algebraic constraint in Eq.~(\ref{eqn:daeProblem}), which ensures that the motion stays on the circle. Parameter $ a $ is a damping coefficient, and $ g $ represents the acceleration due to gravity.   The cost functional in Eq.~(\ref{eqn:daeCostFun}) seeks to minimize the error in tracking a moving target which leads the pendulum in phase by an amount, $\alpha$, with an additional penalty on the control effort. In Eq.~(\ref{eqn:daeCostFun}), parameters $ c$ and $d$ are weights used to emphasize either control-effort or tracking.

The optimal control problem as posed in Eqs.~\eqref{eqn:daeCostFun} and \eqref{eqn:daeProblem} can be solved quite directly using second-order PS differentiation matrices\cite{boyd,atap,RossMED,RossFahroo:issues, RossFahroo:flat}; however, we follow Campbell et al\cite{campbell2016solving} in first transforming the second-order system to a standard state-space form:
\begin{numcases} {(D_I)} \label{eqn:daeFirstOrderSys}
	\dot x_1 = x_2 \nonumber\\
	\dot x_2 = -x_5x_1 - ax_2 + ux_3 \nonumber\\
	\dot x_3 = x_4 \nonumber\\
	\dot x_4 = -x_5x_3 - ax_4 - g - ux_1 \nonumber\\
	\;\;0 = x_1^2 + x_3^2 - L^2  \label{eqn:algebraicConstraint} \\
	x_1(0) = 0 \;\; x_2(0) = L \;\; x_3(0) = 0 \;\; x_4(0) = 0 \nonumber
\end{numcases}	
where variables, $x, \dot x, y, \dot y $ are denoted by $ x_1, x_2, x_3, x_4 $ and $ \lambda $ is now denoted by $ x_5 $.  According to Ref.~[\citen{campbell2016solving}],  the index of the algebraic constraint in Eq.~(\ref{eqn:daeFirstOrderSys}) is  three.

Sophisticated nonlinear programming techniques together with several advanced direct transcription (DT) methods were used in Ref.~[\citen{campbell2016solving}] to solve the DAE optimal control problem.  The success or failure of the various approaches was shown to be strongly dependent on how the problem was transformed for computation.  In particular, Campbell and Kunkel\cite{campbell2016solving} showed that it is crucial to properly transform the problem for DT methods to be successful in the presence of high index DAEs. Various techniques, including index reduction, were applied to problem $ (D) $, some which worked, and others which did not. The reader is encouraged to read [\citen{campbell2016solving}] for an in-depth discussion on the difficulties of solving higher index DAE optimal control processes.

\section{A Pseudospectral Answer to the Challenge Problem}
\begin{figure}[htb]
	\centering
	\makebox[\textwidth][c]{\includegraphics[width=0.91\textwidth]{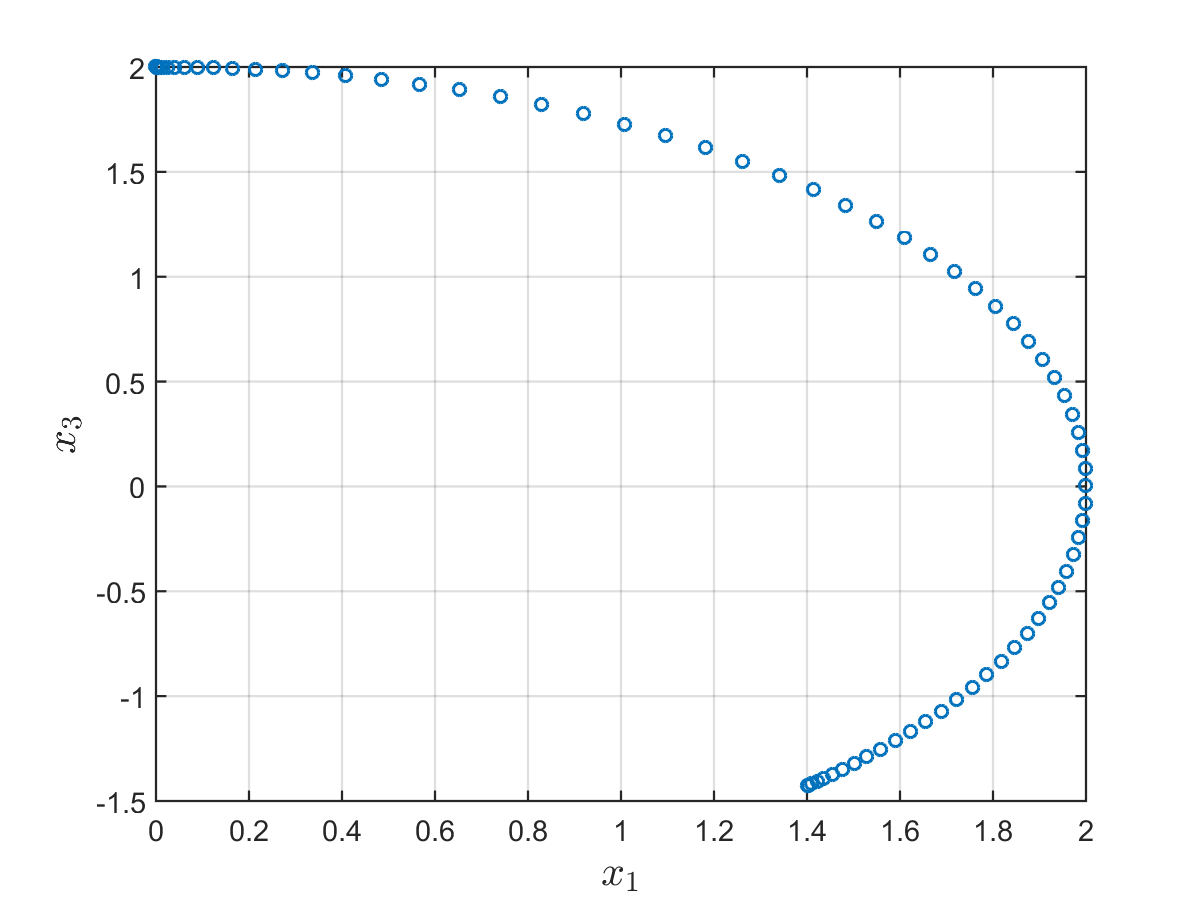}}
	\caption{Candidate optimal state trajectory in $ x_1,x_3 $ plane. }
	\label{fig:xyPhasePlot}
\end{figure}

\begin{figure}[htb]
	\centering
	\makebox[\textwidth][c]{\includegraphics[width=0.91\textwidth]{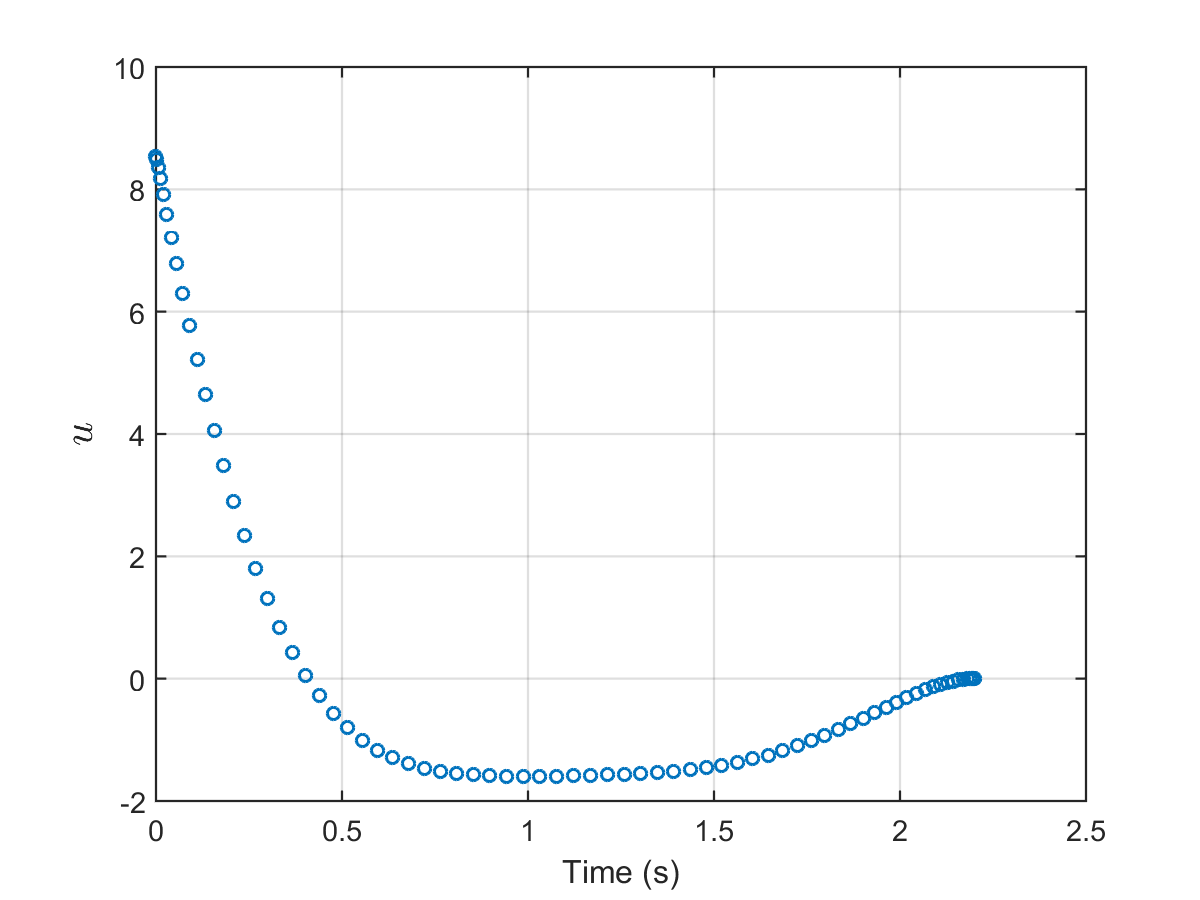}}
	\caption{Candidate optimal control solution. }
	\label{fig:u_tau}
\end{figure}

\begin{figure}[htb]
	\centering
	\makebox[\textwidth][c]{\includegraphics[width=0.9\textwidth]{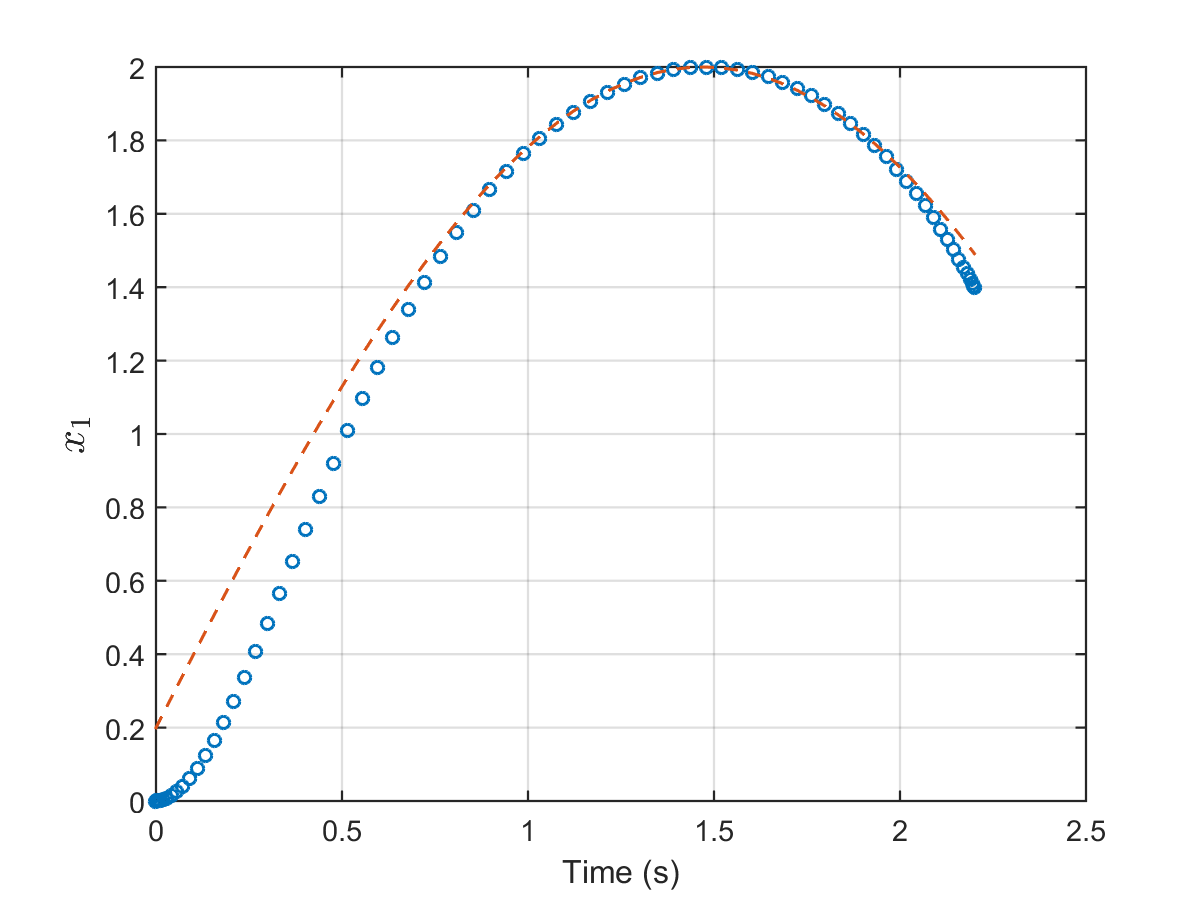}}
	\caption{State solution for $x_1$ with target (dashed line). }
	\label{fig:x1Target}
\end{figure}

\begin{figure}[htb]
	\centering
	\makebox[\textwidth][c]{\includegraphics[width=0.9\textwidth]{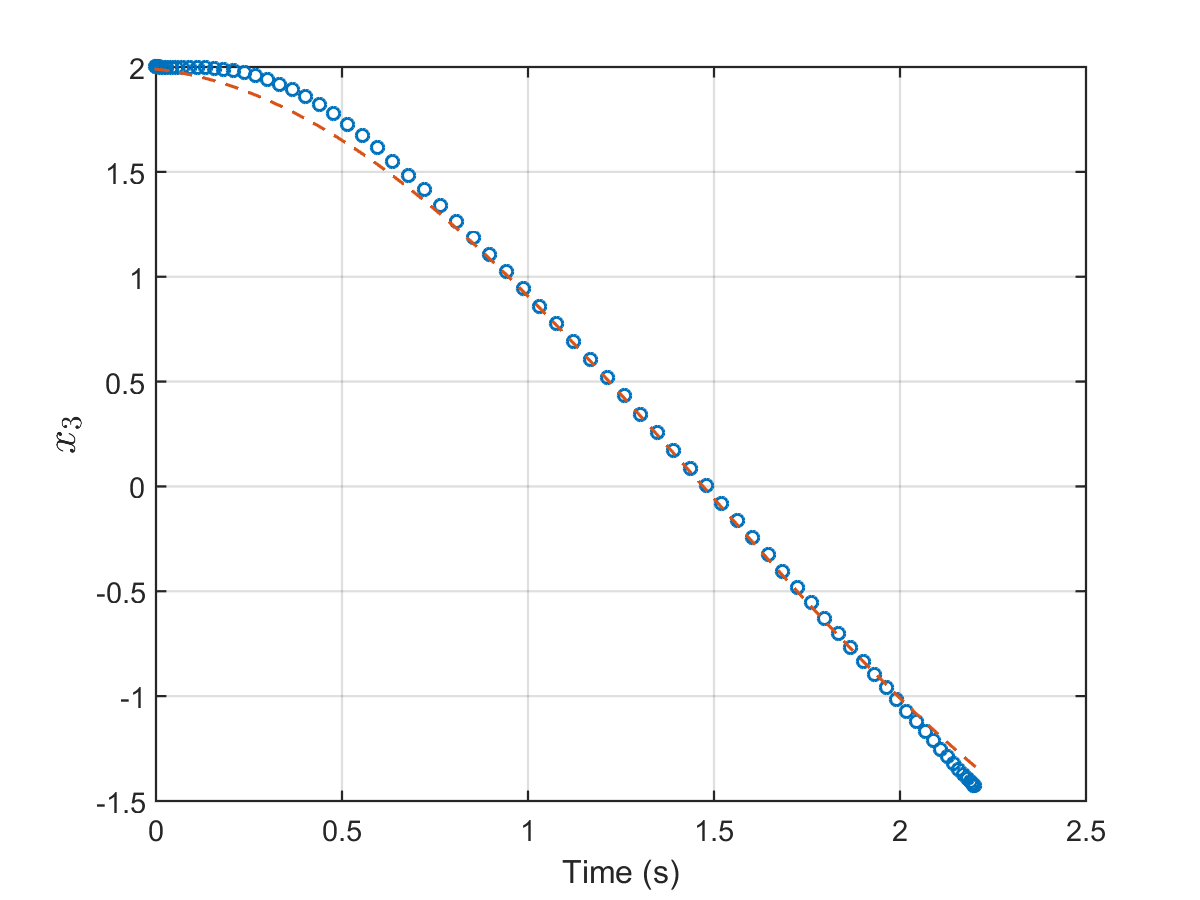}}
	\caption{State solution for $x_2$ with target (dashed line). }
	\label{fig:x2Target}
\end{figure}

The optimal control problem given by Eqs.~(\ref{eqn:daeFirstOrderSys}) and ~(\ref{eqn:daeCostFun}) was coded as given and solved in DIDO. The data  for problem $ (D_I) $ are exactly the same as those used in Ref.~[\citen{campbell2016solving}], and are given as: $ a=0.5, c=1, d=100, g=4, L=2 $, and  $ T=2.2 $. The results are shown in Figures~\ref{fig:xyPhasePlot} through \ref{fig:x2Target}.  A visual inspection of the numerical solution shows that the results obtained using the PS method implemented in DIDO are the same as those presented in~[\citen{campbell2016solving}], where the latter solution was obtained by modifying the original problem given by Eqs.~(\ref{eqn:daeFirstOrderSys}) into a mathematically equivalent form with no path constraints (see Figures 2 and 3 of [\citen{campbell2016solving}]).

It is also worth noting that DIDO does not require a ``guess'' to use the software\cite{scaling}.  In order words, the solution presented in Figures~\ref{fig:xyPhasePlot} through \ref{fig:x2Target} is ``unbiased'' from a user's perspective. Because there was no clear difficulty in generating this solution, our results immediately suggest the need for a new DAE theory of difficulty.

\section{An Independent Verification of Feasibility}
At NASA and elsewhere in the industry, it is not sufficient to generate a solution as presented in Figures~\ref{fig:xyPhasePlot} through \ref{fig:x2Target} without subjecting it to a battery of mathematical and engineering tests.  Furthermore, an argument that the residuals are small (e.g., $10^{-8}$) is well-established to be irrelevant\cite{scaling,RossReview,spec-alg,arb-grid} because it is possible to produce a wrong answer --- even an infeasible one --- with very small ``collocation errors.''  A standard industry practice in testing the accuracy of a control solution is to subject it to an established propagation technique and compute certain critical functions of the propagated state variables; see Refs.~[\citen{ross-book}] and [\citen{scaling}] for further details.  In following this industrial rigor, we interpolate the control solution presented in Figure~\ref{fig:u_tau} to generate $u(t)\ \forall\ t \in [0, T] $ and propagate the initial conditions through the ODE $\dot x = f(x, u(t))$ using \textsf{ode45} in MATLAB.

\begin{figure}[htb]
		\centering
		 \makebox[\textwidth][c]{\includegraphics[width=0.90\textwidth]{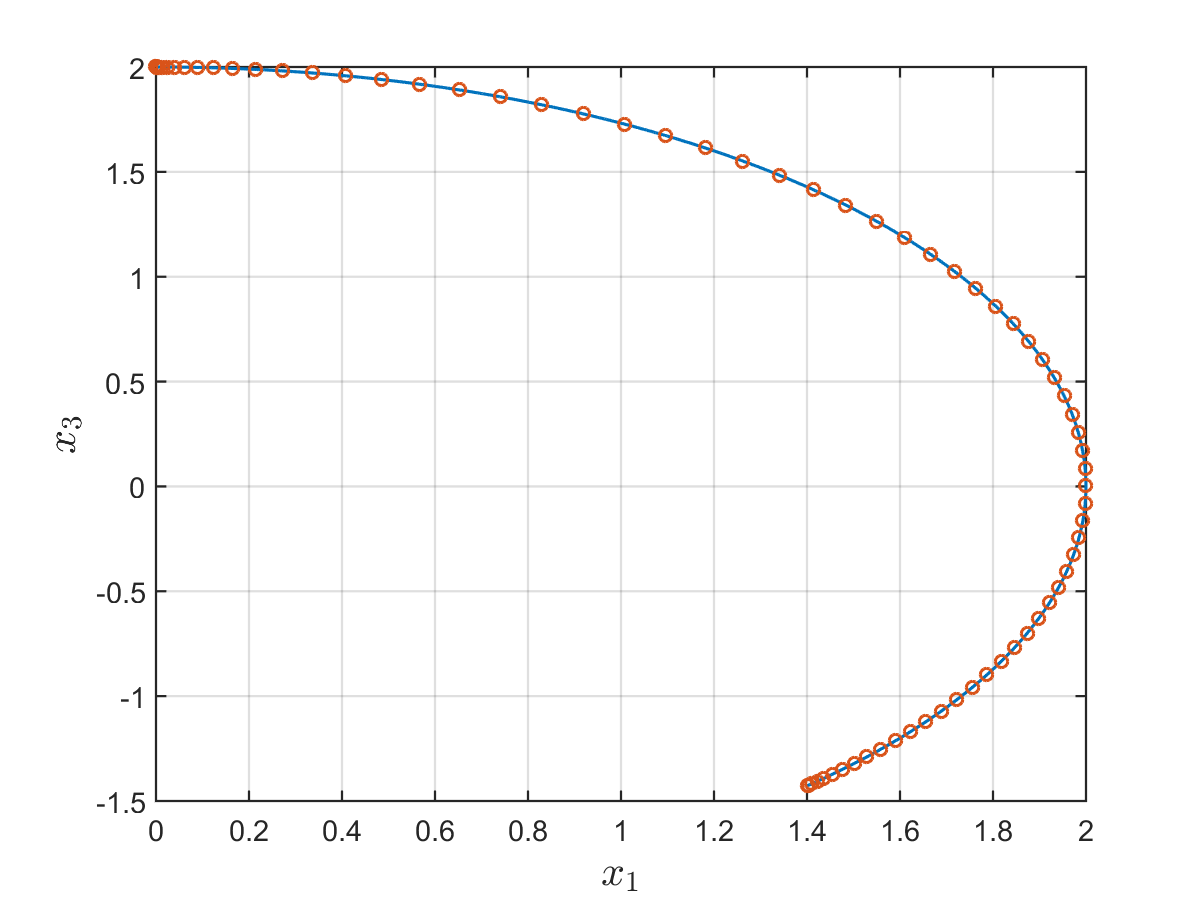}}%
		\caption{Feasibility of a PS solution to problem $ (D_I) $. Legend: data points -- PS solution; solid line -- propagation of the optimal controls.}
		\label{fig:feasState}
	\end{figure}

Figure~\ref{fig:feasState} shows a result from this test.
The DIDO result of Figure~\ref{fig:xyPhasePlot} is overlaid on Figure~\ref{fig:feasState}. It is visually apparent that there is excellent agreement between the propagated solution and the DIDO states.  Numerically, this agreement is within $\pm 1 \times 10^{-4} $.  This number is well within typical industry tolerances for errors; hence we declare the numerical optimal control solution is a dynamically feasible solution to problem $ (D_I) $.

\begin{figure}[htb]
		\centering
		\makebox[\textwidth][c]{\includegraphics[width=0.90\textwidth]{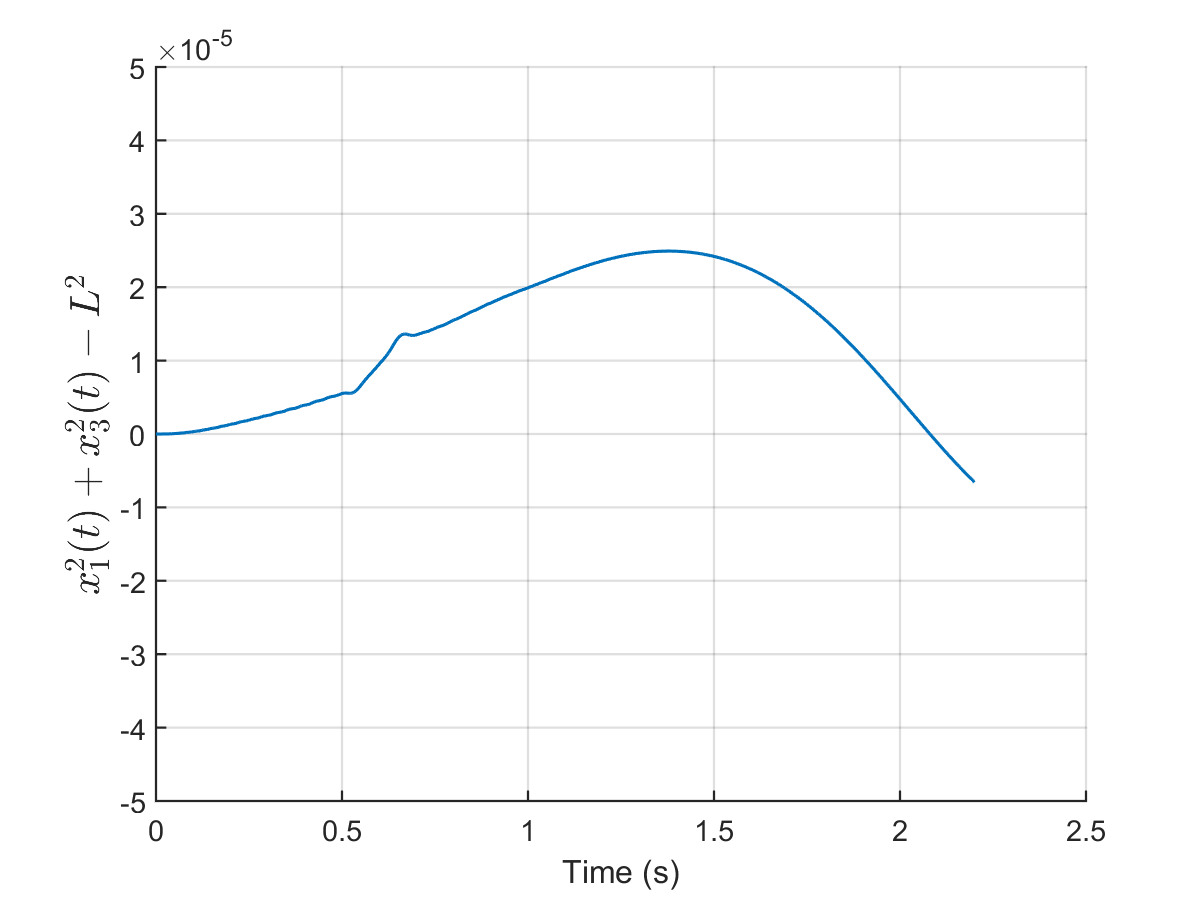}}
		\caption{Satisfaction of algebraic constraint, $ x_1^2+x_3^2-L^2 =0 $, using propagated optimal control.}
		\label{fig:feasPath}
	\end{figure}

Using the propagated values of the state variables, we now evaluate the feasibility of the algebraic constraint, $ x_1^2+x_3^2-L^2 =0 $. Figure~\ref{fig:feasPath} shows that the algebraic constraint, evaluated using the propagated states, is met to well within $\pm 1 \times 10^{-4} $.  As a result, we now declare that that the DIDO solution is an independently-verified, high-quality, DAE-feasible solution to problem $(D_I)$.	

It is important to note that we used RK techniques for propagation; not for optimization. In a typically industry setting, a myriad of alternative verification tools are used to gauge the quality of a solution\cite{Karp-JWST,Kepler-micro-slew,TRACE-IEEE-Spectrum,RossReview}. Furthermore, a determination of the quality of the solution is specific to the application.  For instance, in the case of the Kepler spacecraft, a pointing accuracy of milliarc seconds is mission critical.  This level of accuracy was obtained in [\citen{Kepler-micro-slew}] using a mere 30 PS nodes in DIDO.  In other applications, fewer or higher number of PS nodes may be necessary based on the specific tolerance requirements as gauged by a suite of independent verification tools and not on the nonlinear programming tolerance settings.

\section{Possible Explanations for the Success of PS Tools}
In this section, we offer a suite of possible explanations for the success/failure of PS/RK methods and their various implementations.  We break up our explanations in terms of information that is specific to problem $(D_I)$ and the general differences between PS/RK methods and their software implementations.

\subsection{Necessary Conditions for Optimality}\label{sec:necCondDerivation}
The necessary conditions for optimality of problem $ (D_{I}) $ can be easily obtained by a straightforward application of Pontryagin's Minimum Principle\cite{ross-book}.
At the heart of Pontryagin's Minimum Principle is the Hamiltonian Minimization Condition (HMC).  The HMC states that, for an extremal control $ \bm{u}^* $ to be optimal, it must minimize the control Hamiltonian at each instant of time.  Due to the presence of state-dependent path constraints in problem $ (D_{I}) $, the necessary conditions from the HMC are obtained from the Lagrangian (of the control Hamiltonian):
\begin{equation}\label{eqn:legHamilGeneral}
\overline{H}(\bm{\mu}, \bm{\lambda}, \bm{x}, \bm{u}, t) = H(\bm{\lambda}, \bm{x}, \bm{u}, t) + \bm{\mu}^\intercal \bm{h}( \bm{x}, \bm{u}, t)
\end{equation}
where $ \bm{x}, \bm{\lambda}, \bm{u}, $ represent the vectors of state, costate, and control and are of dimension $ N_x, N_x $ and $ N_u $. The scalar-valued function $ H(\bm{\lambda}, \bm{x}, \bm{u}, t)\in\bbr $ is the control Hamiltonian, $ \bm{\mu}\in\bbr^{N_p} $ the path covectors associated with the HMC, and the vector-function $ \bm{h}( \bm{x}, \bm{u},t)\in\bbr^{N_h} $ corresponds to the path constraint(s) of dimension $ N_h $. A necessary condition is that the Lagrangian of the Hamiltonian  be stationary with respect to the control $ \bm{u} $:
\begin{equation}\label{eqn:statCondGeneral}
	\libzptrl{\overline{H}}{\bm{u}} = \libzptrl{{H}}{\bm{u}} + {\paren{\libzptrl{\bm{h}}{\bm{u}}}}^\intercal\!\bm{\mu}=\bm{0}\in\bbr^{N_u}
	\end{equation}
Along with the stationarity condition, the KKT conditions require the individual path covectors, $ \mu_i $, satisfy the complementarity conditions:
\begin{equation}\label{eqn:compCondGeneral}
	\mu_i  \begin{cases}
	\ge0 \;\; \mbox{if $ h_i(x, u) = h_i^U $} \\
	=0\;\; \mbox{if $h_i^L < h_i(x, u) < h_i^U $} \\
	\le0\;\; \mbox{if $ h_i(x, u) = h_i^L $}
	\end{cases}
\end{equation}
where $  \bm{h}^L,\bm{h}^U\in\bbr^{N_p} $ serve as lower and upper bounds for the vector of path constraints $ \bm{(h)} $, i.e. $ \bm{h}^L\le \bm{h}( \bm{x}, \bm{u},t)\le \bm{h}^U $.
	
Because there is only one path constraint, vector-function $ \bm{h} $ reduces to a scalar function and is further dependent only on the state  so $ h(\bm{x}) = x_1^2+x_3^2-L^2 $. The Lagrangian of the control Hamiltonian is thus given by:
\begin{multline}\label{eqn:legHamil}
	\overline{H}(\bm{\mu}, \bm{\lambda}, \bm{x}, \bm{u}, t) = cu^2 + d\paren{x-L\sin(t+\alpha)}^2 + d\paren{x_3-L\cos(t+\alpha)}^2  \\
	+ \lambda_{x_1}x_2 + \lambda_{x_2}\paren{-x_5x_1-ax_2+u x_3} + \lambda_{x_3}x_4 + \lambda_{x_4}\paren{-x_5x_3-ax_4 - g - u x_1} \\
	+ \mu\paren{x_1^2+x_3^2-L^2}
\end{multline}
where the subscripts on each costate corresponds to the associated state variable.
	
From the stationary condition, the necessary conditions on the controls may be derived. Firstly,
\begin{align}
	\dfrac{\partial \overline{H}}{\partial u} &= u-\dfrac{1}{2c}\paren{\lambda_{x_4}x_1 - \lambda_{x_2}x_3}=0
\end{align}
from which we obtain
	\begin{align}
		\label{eqn:NC01}
		u &= \dfrac{1}{2c}\paren{\lambda_{x_4}x_1 - \lambda_{x_2}x_3}
	\end{align}
	Secondly, we have
	\begin{align}
		\label{eqn:NC02}
		\dfrac{\partial \overline{H}}{\partial x_5} &= -\lambda_{x_2}x_1 -\lambda_{x_4}x_3=0
	\end{align}
	Because $ x_5 $ appears linearly in Eq.~(\ref{eqn:legHamil}), the optimal control is singular\cite{bryson1975applied}.  Nonetheless in the next section, Eq.~\ref{eqn:NC02} will be used to verify the optimality of $x_5$.

Due to the presence of the path constraint the adjoint variables for problem $ (D_I) $ evolve according to $ -\dot{\bm{\lambda}}  = \dfrac{\partial \overline{H}}{\partial \bm{x}}$. The adjoint equations are given as follows:
	\begin{align}
	-\dot\lambda_{x_1} &= 2d\paren{x_1-L\sin(t+\alpha)}-\lambda_{x_2}x_5-\lambda_{x_4}u+2\mu x_1 \label{eqn:costateFirst}\\
	-\dot\lambda_{x_2} &= \lambda_{x_1}-a\lambda_{x_2} \\
	-\dot\lambda_{x_3} &=  2d\paren{x_3-L\cos(t+\alpha)}+\lambda_{x_2}u -\lambda_{x_4}x_5 + 2\mu x_3\\
	-\dot\lambda_{x_4} &= \lambda_{x_3}-a\lambda_{x_4} \label{eqn:costateLast}
	\end{align}
	Finally, the terminal transversality conditions are given by $ \bm{\lambda}(t_f)=  \libzptrl{\overline{E}}{\bm{x}(t_f)}$.  Since problem $ (D_I) $ only specifies initial values, the terminal state does not appear in the endpoint Lagrangian. Hence, the value of all the costates must be zero at the final time $ T $, i.e.
	\begin{equation}\label{eqn:terminalTransversalityCond}
	\left [ \lambda_{x_1}(T), \lambda_{x_2}(T), \lambda_{x_3}(T), \lambda_{x_4}(T) \right ] = \mathbf{0}
	\end{equation}	

\subsection{Verification of the Necessary Conditions via DIDO}
In addition to the primal variables, DIDO also outputs a suite of dual variables.  This does not imply that DIDO implements and ``indirect'' PS method.  Its interface to the user is only the ``direct'' problem formulation.  Because the necessary conditions are eventually procedural, DIDO ``derives'' these conditions and implements it via the covector mapping principle (CMP)\cite{ross-book}.  Hence, the recommended approach\cite{ross-book} to using DIDO is to apply it a ``mathematical tool'' for solving optimal control problems.  In other words, if the problem formulation or its necessary conditions exhibit certain pathological behavior, DIDO will reflect it.  Typical pathologies to avoid in using DIDO are unboundedness (e.g., variables going to infinity inside the search space) and nonsmoothness in data or its Jacobians (e.g., the square root function). \emph{\textbf{Singular arcs and DAEs are not considered pathological}}. Hence, from an analysis of the necessary conditions derived in Section~\ref{sec:necCondDerivation}, it is apparent that DIDO should work.

\begin{figure}[htb]
		\centering
		\makebox[\textwidth][c]{\includegraphics[width=0.9\textwidth]{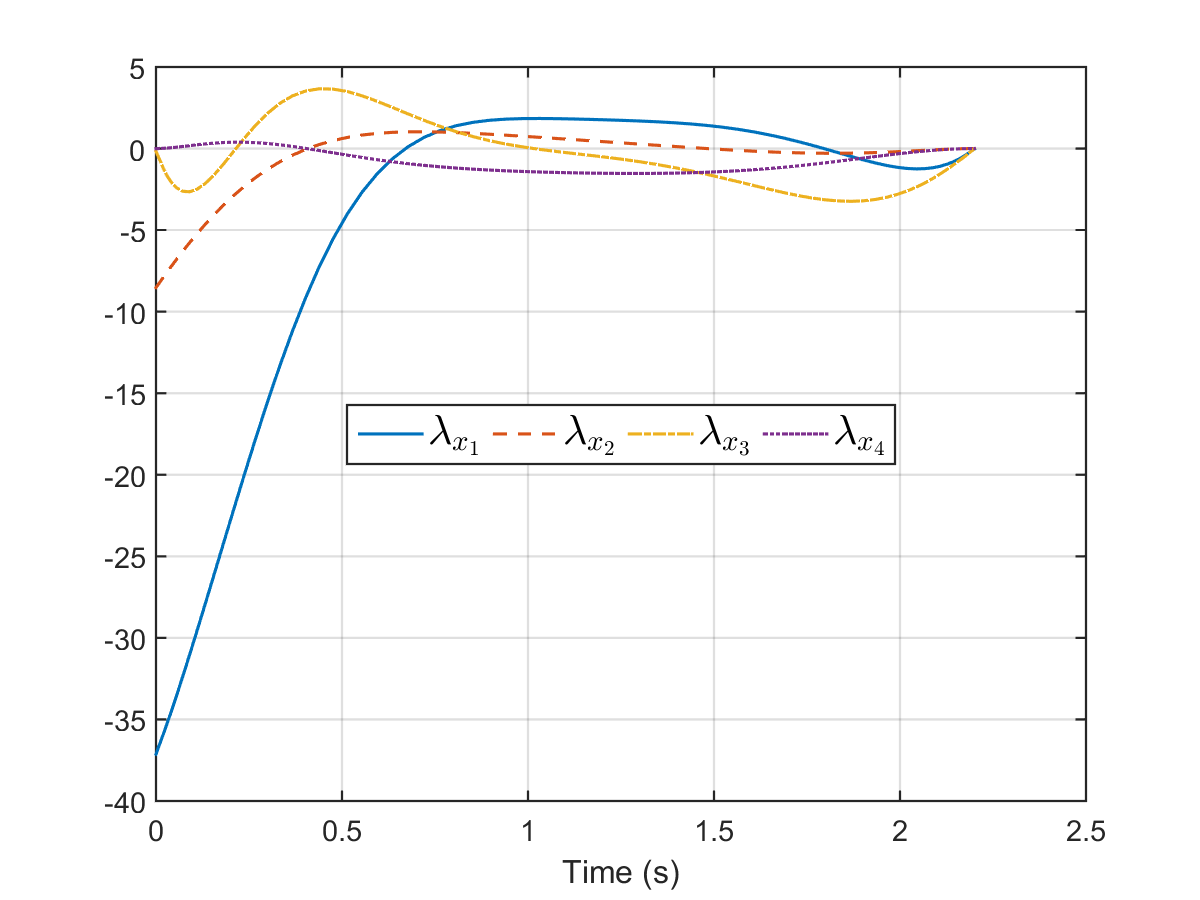}}
		\caption{Costates associated with the state variables. }
		\label{fig:costates}
\end{figure}

A plot of the DIDO-generated costates is shown in Figure~\ref{fig:costates}. It is clear from this figure that all the costates take the value zero at the terminal time. In other words, the terminal transversality condition given by Eq.~\eqref{eqn:terminalTransversalityCond} is satisfied numerically.

\begin{figure}[htb]
		\centering
		 \makebox[\textwidth][c]{\includegraphics[width=0.90\textwidth]{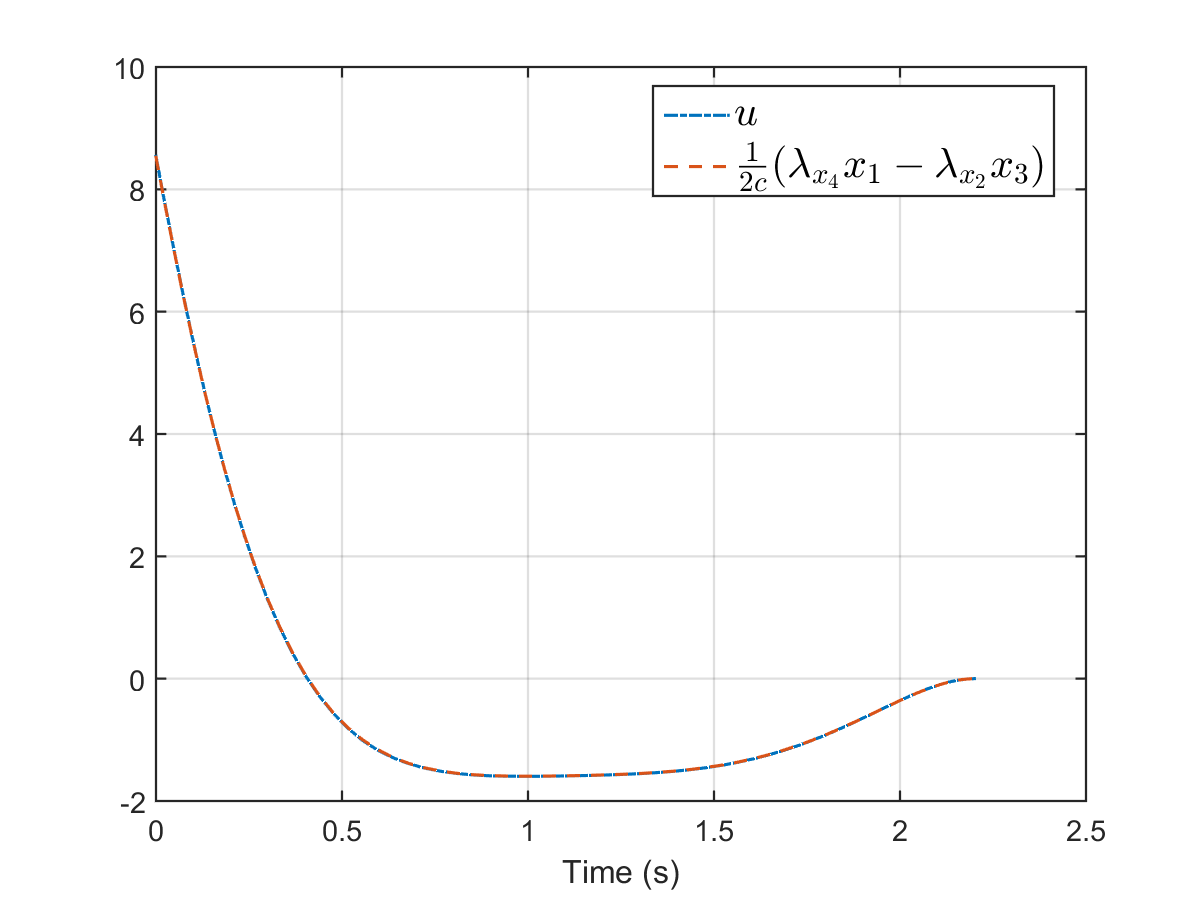}}
		\caption{Verification of the necessary condition in Eq.~(\ref{eqn:NC01}).}
		\label{fig:statCond_ctrlTau}
	\end{figure}

\begin{figure}[htb]
		\centering
		\makebox[\textwidth][c]{\includegraphics[width=0.90\textwidth]{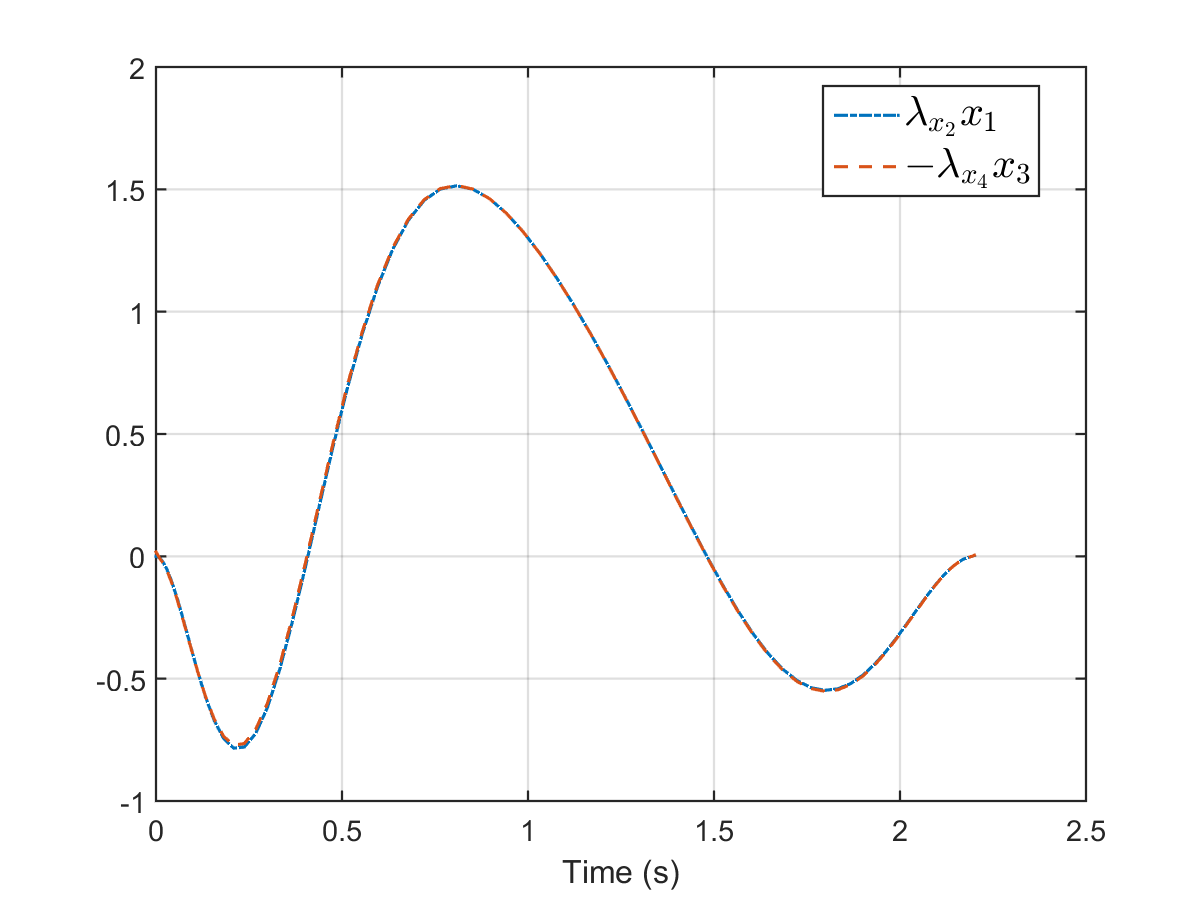}}
		\caption{Verification of the necessary condition in Eq.~(\ref{eqn:NC02}).}
		\label{fig:statCond_ctrlF}
	\end{figure}

Equation~(\ref{eqn:NC01}) gives an expression for the control, $u$, in terms of the costates. Plotting both sides of the equation, as is done in Figure~\ref{fig:statCond_ctrlTau}, illustrates that the numerical solution adheres to the necessary condition on the regular control. Similarly, we plot $\lambda_{x_2}x_1$ against $-\lambda_{x_4}x_3$ to check the condition on the singular control given by Eq.~(\ref{eqn:NC02}). Figure~\ref{fig:statCond_ctrlF} shows that Eq.~(\ref{eqn:NC02}) is also satisfied by the numerical solution to problem $D_I$. Thus, we may reasonably conclude that the control presented in Figure~\ref{fig:u_tau} is an extremal solution to problem $D_I$.  Furthermore, we note that there was no major issue in solving a DAE problem with a singular arc.

\subsection{An Overview of Why PS Theory Works}
Since the year 2007, PS optimal control methods have rapidly evolved to the status of a fundamentally new approach to optimal control theory itself. See Ref.~[\citen{RossReview}] and the references contained therein for details of this concept.  In broad terms, PS optimal control theory is based on two fundamentals\cite{RossReview}: (\emph{i}) a continuous state-trajectory can be approximated to any precision by a sufficiently high-order polynomial, and (\emph{ii}) the covector trajectories associated with an optimal control problem can be approximated to spectral accuracy by a transformation of the multipliers associated with the discrete mathematical programming problem.  The first point is a direct consequence of the Stone-Weierstrass theorem\cite{royden1988real,RossReview} , and the second point has been established as part of the CMP\cite{RossReview,gong2008connections}.  Furthermore, convergence theory of PS optimal control does not rely on the CMP, rather it is based on a combination of Polak's notion of consistent approximation\cite{polak2012optimization} and the classic Arzel\`a-–Ascoli theorem~\cite{rudin1964principles}.  The reader is directed to~[\citen{RossReview,TAC:linearizable,kang-ijrc,kang-rate,Kang_2008_convergence}] for details.
Singular arcs and the index of a DAE do not directly enter in the proof of convergence of PS methods.  Consequently, such issues do not appear to be detrimental to a proper implementation of PS methods.  In contrast, proofs of convergence of RK methods rely on very strong assumptions that are frequently not satisfied in many practical applications.
	
\section{Conclusion} \label{sec:conclusion}
Great care must be exercised in drawing \emph{\textbf{negative}} conclusions from software tools and/or na\"{\i}ve implementations of computational methods.  For instance, it is very easy to implement a ``bad'' RK or PS method by deliberately or inadvertently choosing incorrect/inaccurate coefficients, grids, weights etc. Even if sophisticated nonlinear programming solvers are ``patched'' to inappropriate discretization techniques, the resulting ``advanced method'' can be shown to fail for the simplest of the problems. Hence, if a particular implementation does not work, it is inappropriate to conclude a negative result unless it is consistent with theoretical predictions.  In the same spirit, if a particular approach routinely provides consistent results that ``defies theory,'' then it is the theory that must be questioned.

Singular arcs and DAEs are frequently encountered in practical industry applications. These problems have been routinely solved over the last two decades using proper implementations of PS optimal control techniques.  Because of the higher failure rates of RK methods, PS methods have gradually replaced legacy optimization methods, particularly in new and emerging problems in aerospace engineering. It is quite possible that one may be able to construct an innocuous DAE problem that cannot be solved by any PS method.  Consequently, a new theory of difficulty for DAEs is warranted.

\end{document}